\documentclass[aps,twocolumn,prb,floatfix,amsmath,amssymb,groupedaddress]{revtex4}

\usepackage{graphicx}
\usepackage{psfrag}
\usepackage{epstopdf}
\usepackage{bm}
\usepackage[justification=raggedright]{caption}
\usepackage{subfig}
\usepackage{enumitem}
\usepackage[export]{adjustbox}
\usepackage{mdwlist}
\usepackage{amsthm}
\usepackage{amsmath}
\usepackage{multirow}
\newtheorem{prop}{Proposition}

\makeatletter
\RequirePackage{keyval}

\begin{document}

\providecommand{\half}{{\frac{1}{2}}}		

\title{``Irregularization'' of Systems of Conservation Laws}

\author{Hunter Swan}
\author{Woosong Choi}
\author{Matthew Bierbaum}
\author{Yong S. Chen}
\author{James P. Sethna}
\affiliation{Laboratory of Atomic and Solid State Physics (LASSP),
Clark Hall, Cornell University, Ithaca, New York 14853-2501, USA}
\author{Stefanos Papanikolaou}
\affiliation{Department of Mechanical Engineering,
The Johns Hopkins University, Hopkins Extreme Materials Institute,
113 Malone Baltimore, MD, 21218}

\begin{abstract}
We explore new ways of regulating defect behavior in systems of conservation laws.  Contrary to usual regularization schemes (such as a vanishing viscosity limit), which attempt to control defects by making them smoother, our schemes result in defects which are \textit{more singular}, and we thus refer to such schemes as ``irregularizations''.  In particular, we seek to produce \textit{delta shock} defects which satisfy a condition of \textit{stationarity}.  We are motivated to pursue such exotic defects by a physical example arising from dislocation dynamics in materials physics, which we describe.
\end{abstract}

\maketitle

A remarkable feature of systems of conservation laws is the tendency of smooth initial data to give way to non-smooth defects after only finite time.  Examples include shocks in hydrodynamics, cracks and dislocations in crystalline solids, and traffic jams in continuum models of traffic flow.  Such defects are important in determining the qualitative behavior of a system, but they pose various theoretical and practical challenges, stemming from the fact that the time evolution of defects cannot be determined solely from the evolution of the ambient continuum.


The physical origin of this difficulty can be understood as follows:  Continuum equations of motion are derived from more fundamental, microscopic laws by discarding the ``high-order'' terms which are irrelevant at macroscopic scales.  Thus, systems with different microphysics can yield the same continuum equations, and conversely the continuum equations alone do not specify the microphysics.  Defects are inherently microscopic, and thus their behavior cannot be inferred from the continuum laws.  

The ambiguity of defect behavior in systems of conservation laws is reflected in the mathematical theory by \textit{non-uniqueness} of solutions containing defects.  To obtain uniqueness, one must impose additional conditions upon the solution, presumably corresponding to the microphysics not accounted for by the continuum equations.  Common conditions of this sort include \textit{entropy conditions} and any of various types of \textit{regularization}.  In the former case, inequalities are imposed in analogy with the thermodynamical principle that entropy must always increase, while in the latter case the conservation laws themselves are modified by the addition of small terms that serve to smooth out defects.  We will return in more detail to both of these cases in the next section. 

For the most common kinds of conservation laws (notably, those of hydrodynamics), the ``proper'' choice for solutions with defects turns out to be relatively tame, with defects consisting only of jump discontinuities and being described by regularization with a vanishingly small viscosity term.  Some more exotic systems feature different classes of defects which require other types of regularization; for example, superfluids are known\cite{Meppelink2009} to exhibit \textit{dispersive shock waves} which are described by a very small dispersive term.  Still another type of defect (which is a primary subject of this paper) is a \textit{delta shock}, which consists of a delta function-like spike, possibly riding atop a regular (viscous) shock.  This latter type of defect is much more poorly understood than the former two and differs in some important regards: Firstly, there is no particular regularization associated with delta shocks (for some systems, delta shocks may even be contained in the vanishing viscosity solutions\cite{Tan1994}).  Secondly, there are few physically-realizable systems which are currently known to exhibit delta shocks.

\begin{figure}[h]
	\includegraphics[width=1.14\columnwidth,center]{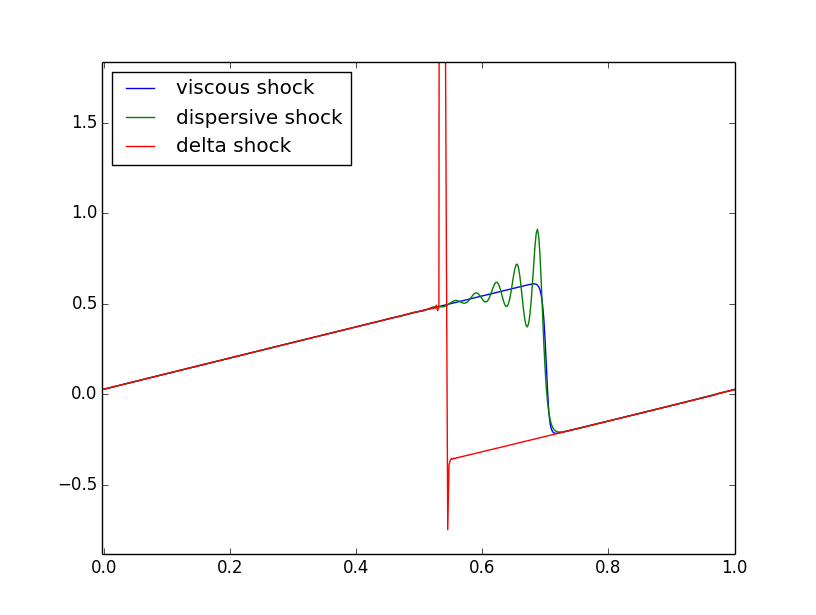}
	\caption{{\bf Shocks in Burgers' equation.} (Color online) Various types of shocks evolving from smooth, sinusoidal initial conditions, depending on different regularizations of the Hopf equation: viscous shock (blue), dispersive shock (green), and delta shock (red). 
	}
	\label{fig-shocks}
\end{figure}

The present paper was motivated by the discovery of delta shocks in a system arising from materials physics---continuum dislocation dynamics (CDD)---which to our knowledge provides the best physical interpretation of delta shocks currently available. These delta shocks have presented substantial challenges to physicists interested in simulating CDD numerically, which has led to a desire for new numerical methods tailored to better handle such exotic defects.  We explore in this paper (in a mostly informal way) perspectives on delta shocks which may be of use to practitioners interested in simulating them numerically.  (We do not, however, explicitly construct any numerical schemes here.)  In particular, we describe a condition of \textit{stationarity} which forces the formation of delta shocks, and we construct a new kind of ``irregularization'' which also produces delta-like defects (though we do not rigorously prove these to be true delta shocks).

The outline of this paper is as follows:  In the first section we provide some background on classical shocks.  In the second section we describe delta shocks and discuss various ways in which they can arise in systems of conservation laws, and we introduce the condition of stationarity and the irregularization mentioned above.  Then in the third section we describe in more detail the system which initiated the preceding work---continuum dislocation dynamics.

\section{Background}
In this section, we review basic facts about systems of conservation laws and standardize our notation.  A system of conservation laws (or simply ``conservation law'') is a partial differential equation (PDE) of the form

\begin{equation}
u_t + \nabla \cdot F(u) = 0,
\label{SCL}
\end{equation}
where $u:\Omega \times \mathbb{R}^+ \rightarrow \mathbb{R}^m$ (for $\Omega\subset \mathbb{R}^n$), $u=u(x,t)$, is the unknown, and $F:\mathbb{R}^m\rightarrow\mathbb{R}^n$ is a given \textit{flux function}.  The conserved quantity $u$ we refer to as ``mass'' in the following discussion.  In the case where $u$ is scalar valued ($m=1$), we refer to eq.~(\ref{SCL}) as a scalar conservation law.  

Our working example of a conservation law will be the Hopf equation (a.k.a. the ``inviscid Burgers' equation''), 
\begin{equation}
u_t + \left(\frac{1}{2} u^2\right)_x = u_t + u u_x = 0,
\label{Hopf}
\end{equation}
which is a scalar conservation law in one spatial dimension ($n=1$).  The Hopf equation is one of the simplest examples of a PDE which develops discontinuities in finite time\cite{Evans}.  Such discontinuities and other non-smooth features appearing in solutions to PDE we refer to collectively as \textit{defects}.  In the presence of discontinuities in $u$, the derivatives appearing in eqs. (\ref{SCL}-\ref{Hopf}) are no longer well defined, and so we must generalize our notion of ``solution'' to a PDE in order to make sense of the conservation law at later times.  Such a generalization is referred to as a \textit{weak solution}.  Several distinct notions of weak solution exist, such as integral solutions, viscosity solutions, variational solutions, and various types of regularized solutions, to name a few\cite{Evans}.  

Not all forms of weak solution are applicable to all PDE's.  The types of weak solution which are most relevant for systems of conservation laws are \textit{regularized solutions} and \textit{integral solutions}.  In the former case, the PDE is modified by the addition of a small term that serves to smooth out the defects; for the Hopf equation, such a modification looks like
\begin{equation}
u_t + \left(\frac{1}{2}u^2\right)_x = \eta \; \partial_x^k u,
\label{RegHopf}
\end{equation}
where $\eta\in\mathbb{R}$ controls the strength of the regularization.  To recover a solution to the original system of conservation laws, we take a limit as $\eta \rightarrow 0$.  Note that, in principle, the form of the regularizing term should correspond to the microphysics of a given system---so for example, systems with viscosity have a second-order regularization, $\eta \partial_x^2 u$, whereas inviscid systems (like plasmas or superfluids) typically have dispersive (third-order) regularization, $\eta \partial_x^3 u$.  Below, in eqs. (\ref{modHopf1}-\ref{modHopf3}), we will consider a different form of regularization altogether, which nonetheless bears a familial resemblance to eq. (\ref{RegHopf}). 

The main alternative kind of weak solution, the integral solution, proceeds instead by integrating the original PDE to remove the problematic derivatives.  Integral solutions are not of direct interest to us here, but we point out one tangential relationship to our current discussion: Integral solutions are generally not \textit{unique}, and to select a proper weak solution requires additional microphysical information.  This information is usually specified as a relationship between the values of the solution $u$ on either side of a defect; when the relationship takes the form of an inequality it is called an \textit{entropy condition}\cite{Lax1973,Lax1957}, and when it is an equality it is called a \textit{kinetic relation}\cite{LeFloch2002}.  The stationarity condition we propose below is similar to a kinetic relation, but with the generalization that the \textit{speed} of the defect can also figure into the equality. 

\section{Delta shocks}
For our purposes, a delta shock (also written ``$\delta$-shock'') may be defined informally as any defect which has a finite amount of mass concentrated at a point.  We note in passing that there are also meaningful notions of $\delta '$-shocks (and $\delta^{(n)}$-shocks, more generally) corresponding to defects containing \textit{derivatives} of delta functions\cite{Panov2006}. Our numerical 
solutions (Section~\ref{subsec:Regularizations}, Figs~\ref{Hopf}
and~\ref{coshFig})
appear to contain an upward and downward-pointing singularity -- the sum
being the strength of the delta shock and the difference a $\delta'$
component. In this section we shall derive conditions for the evolution
of the strength of the delta shock component.

Delta shocks are known to arise naturally in a number of systems of conservation laws---for example, the system 
\begin{eqnarray}
u_t + (u^2)_x = 0 \nonumber \\
v_t + (uv)_x = 0
\label{TansEq}
\end{eqnarray}
considered by Tan, Zheng, \& Zhang\cite{Tan1994}, and the system
\begin{eqnarray}
u_t + (u^2-v)_x = 0 \nonumber \\
v_t + \left(\frac{1}{3}u^3-u\right)_x = 0
\label{KeyfitzEq}
\end{eqnarray}
studied by Keyfitz \& Kranzer\cite{Keyfitz1989}.  In the former case, the delta shocks have been shown to be vanishing viscosity limits of the regularized system
\begin{eqnarray}
u_t + (u^2)_x = \eta \; u_{xx} \nonumber \\
v_t + (uv)_x = 0 .
\label{TansViscEq}
\end{eqnarray}

The mechanisms which cause delta shocks to form in systems of conservation laws are not fully understood\cite{Keyfitz1999}, though in some cases the necessity of delta shocks is obvious.  We highlight two such cases here: Firstly, if a conserved quantity $v$ is a derivative of another variable $u$, and the latter exhibits regular (viscous) shocks, then clearly $v$ must contain a delta shock.  For example, if we take a derivative of the Hopf equation (\ref{Hopf}) and set $v:=u_x$, we find that
\begin{equation}
v_t + (uv)_x = u_{xt} + (uu_x)_x = \partial_x \left(u_t + \left(\frac{1}{2}u^2\right)_x \right) = 0,
\label{Hopfddx}
\end{equation}
which is precisely the latter half of system (\ref{TansEq}).  Thus (provided the initial data for system (\ref{TansEq}) satisfies $v(x,t=0)=u'(x,t=0)$) we see that this system can be interpreted as the evolution of the conserved quantity of the Hopf equation, along with that of its derivative.  

This interpretation is not the only available, nor is it necessarily the best.  Indeed, the closely related system\cite{Tan1994}
\begin{eqnarray}
u_t + (F(u))_x = 0 \nonumber \\
v_t + (g(u)v)_x = 0,
\label{modTansEq}
\end{eqnarray}
where $F$ is an arbitrary smooth \& convex function and g is an arbitrary smooth \& increasing function, also displays delta shocks, though they cannot in this case be interpreted as derivatives of regular shocks.  A more general interpretation for the appearance of delta shocks in such conservation laws is as a consequence of \textit{incompatibility of defect motion and the flow of the conserved quantity}, which we explore in the next two subsections. 

\subsection{Generalized Rankine-Hugoniot relation}
The Rankine-Hugoniot relation is a statement of conservation of mass across a shock.  For the case of a regular shock in one spatial dimension which moves with speed $\sigma$, mass conservation requires that the net flux into the shock balance the mass of the ``hole'' left as the shock moves away.  In symbols, this reads
\begin{equation}
F(u_l)-F(u_r) = \sigma \left(u_l - u_r\right),
\label{RH}
\end{equation}
where $u_l$ and $u_r$ are the values of $u$ to the immediate left and right (resp.) of the shock (see Evans\cite{Evans} for a more formal derivation of this relation).

If we allow the possibility of delta shocks storing mass on the defect, it is now the combination of the hole behind the shock \textit{and} the mass $m$ of the delta which must balance the incoming flux.  The modified Rankine-Hugoniot relation is thus
\begin{equation}
F(u_l)-F(u_r) = \sigma \left(u_l - u_r\right) + \frac{dm}{dt}.
\label{MRH}
\end{equation}

Note that both eqs. (\ref{RH}) and (\ref{MRH}) are vector identities in general, and may be interpreted as specifying the shock speed $\sigma$ and delta mass $m$ (or the derivative thereof) in terms of the local environment, given by $u_l$ and $u_r$.  In particular, this implies that eq.~(\ref{RH}) is only satisfiable in general when $u$ is scalar, since vector $u$ leads to an overspecification in (\ref{RH}) of the single unknown $\sigma$.  Eq.~(\ref{MRH}) has no such problem, since $m$ has the same number of components as $u$, so that there is always one more unknown than there are equations in (\ref{MRH}).  This shows that, except in the presence of high degeneracy so that eq.~(\ref{RH}) is soluble for $\sigma$, we expect multicomponent conservation laws to require delta shocks in order to conserve mass in the vicinity of a defect. 

This provides our preferred explanation of the delta shocks in system (\ref{TansEq}): the classical Rankine-Hugoniot relation (\ref{RH}) for this system reads:
\begin{eqnarray}
u_l^2-u_r^2 = \sigma \left(u_l - u_r\right) \implies \sigma = u_l+u_r \nonumber \\
u_l v_l - u_r v_r = \sigma \left(v_l - v_r\right) \implies \sigma = \frac{u_l v_l - u_r v_r}{v_l - v_r},
\end{eqnarray}
which does not hold unless $u_l v_r = u_r v_l$ (i.e. the necessary degeneracy for (\ref{RH}) to be satisfiable does not generally exist).  A similar statement shows why the more general system (\ref{modTansEq}) also displays delta shocks. 

\subsection{Stationarity}
\label{sec:Stationarity}
The Rankine-Hugoniot relations illustrate how multicomponent systems can develop delta shocks when a defect's motion cannot be simultaneously compatible with the flow of all components of the conserved quantity.  We can also produce a similar effect even in scalar conservation laws by insisting that a defect move according to specified kinematics which are not compatible with the flow of the scalar conserved quantity.  The simplest such example is to demand that a defect be \textit{stationary}.  This might be an appropriate kinematic constraint in, for example, an electrical network containing fuses, where a blown fuse is modeled as a stationary defect.  A stationarity condition could also be physically appropriate for the continuum dislocation dynamics example we will describe below.  This stationarity condition is analogous to a kinetic relation\cite{LeFloch2002}, with the generalization that the speed of the defect can now figure in as well (rather than just the value of the local conserved quantity).  

As an illustration, we consider the Riemann problem 
\begin{equation}
u(x,t=0) = \left\{ 
\begin{array}{l l}
1 & \quad x\leq 0\\
0 & \quad x>0 
\end{array} 
\right.
\end{equation}
for the Hopf equation (\ref{Hopf}).  Mass flows into the defect at $x=0$ with rate $F(0^-) - F(0^+) = \frac{1}{2}\cdot1^2-\frac{1}{2}\cdot0^2 = \frac{1}{2}$, and thus by the modified Rankine-Hugoniot relation (\ref{MRH}) with $\sigma=0$ we find that the mass of the delta at $x=0$ grows at a rate $\dot{m} = \frac{1}{2}$.  Everywhere else the value of $u$ is uniquely specified by projection of characteristics.  Hence a reasonable candidate solution for all times $t \ge 0$ is 
\begin{equation}
u(x,t) = u(x,t=0) + \frac{t}{2}\:\delta(x)
\label{HopfSol}
\end{equation}

We note that this candidate solution is not an integral solution in the usual sense, as the square of a delta function (needed to make sense of an integral solution to eq.~(\ref{Hopf})) is undefined by classical distribution theory.  (However, it is possible to overcome this technicality using generalized distributions\cite{Colombeau1984}.)

\subsection{Regularizations for stationarity}
\label{subsec:Regularizations}
Our \textit{ad hoc} construction of a solution (\ref{HopfSol}) to the Hopf equation which satisfies the condition of stationarity does not lend itself well to numerical simulation.  We thus find it advantageous to try to construct a ``irregularization'' of the Hopf equation which yields such a solution.  Our attempted irregularization looks like:
\begin{equation}
u_t + \left(\frac{1}{2}u^2 f(\epsilon,u_x)\right)_x= 0,
\label{modHopf1}
\end{equation}
where the smooth \& bounded function $f(\epsilon,u_x)$ satisfies:
\begin{enumerate*}
\item $f(\epsilon,0) = 1$ \\
\item $f(\epsilon,-\infty) = 0$ \\
\item $f(0,\cdot ) = 1$.
\end{enumerate*}

For example, 
\begin{equation}
u_t + \left(\frac{\frac{1}{2}u^2}{1+\epsilon \: u_x^2}\right)_x = 0
\label{modHopf2}
\end{equation}
is such an irregularization.
The intuition behind this is that when a shock develops in the Hopf equation, $u_x$ diverges to $-\infty$, sending the modified flux $\frac{1}{2}u^2 f(\epsilon,u_x)$ to zero.  The quenching of the flux will hopefully cause a pileup of mass behind the defect.  On the other hand, for small values of $\epsilon$ and away from defects, condition 3 above implies that $f\approx 1$, and the original Hopf equation is recovered (approximately).  

(This prescription for the irregularization is not as general as possible---for example, one also could add to the denominator of eq.~(\ref{modHopf2}) a non-vanishing polynomial in $u$, which would reflect some sort of self-interaction effect.  However, the ``minimal'' irregularization of eq.~(\ref{modHopf2}) is sufficient to produce the defects we are interested in.  Moreover, the effect is not peculiar to this particular form: We have also experimented with the irregularization
\[u_t + \left(\frac{1}{2}u^2 \exp\left(-\epsilon \: u_x^2\right)\right)_x = 0\]
and seen the same qualitative behavior.
)

For numerical work, it is sometimes helpful to also add a small regularizing term to ensure solutions stay relatively smooth---e.g. 
\begin{equation}
u_t + \left(\frac{\frac{1}{2}u^2}{1+\epsilon \: u_x^2}\right)_x = \eta \: \partial_x^k u,
\label{modHopf3}
\end{equation}
where k is typically 2 or 4.

\begin{figure}[h]
	\includegraphics[width=1.14\columnwidth,center]{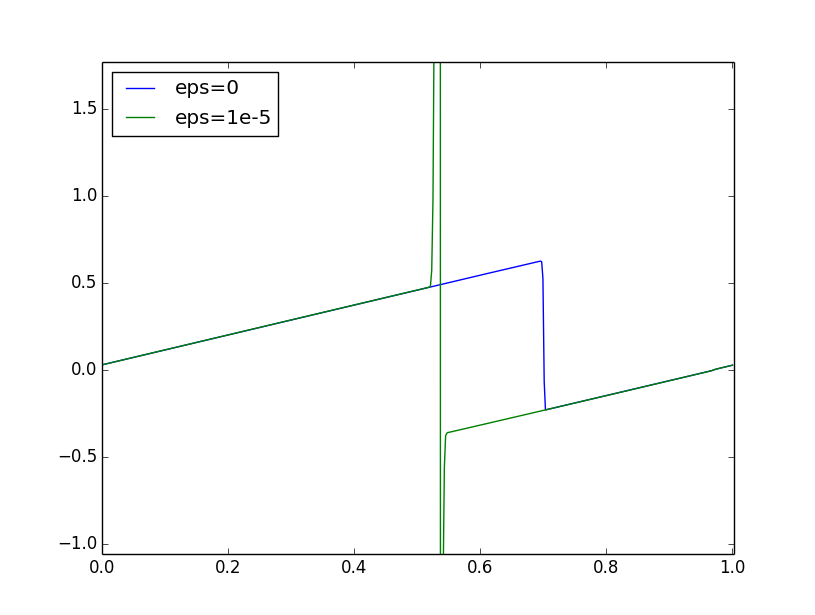}
	\caption{{\bf Irregularized Hopf equation.} (Color online) Simulation of the irregularized Hopf equation (\ref{modHopf3}) starting with sinusoidal
initial conditions, for two values of $\epsilon$ at time $t=1$ (simulated with a conservative upwind differencing). Note that the singularity has both a positive and negative peak. Only the sum contributes to the delta shock (as a weak limit); the difference can be described as a derivative of a $\delta$-function.
	}
	\label{fig-Hopf}
\end{figure}

Does this succeed in making the defect stationary?  We begin to answer this by examining numerical results.  Fig.~(\ref{fig-Hopf}) shows simulations of eq.~(\ref{modHopf3}) for $\epsilon=0$ \& $\epsilon=10^{-5}$ on the unit interval with periodic boundary conditions, for parameters $t=1$, $k=2$, $\eta=10^{-6}$, and a grid spacing $dx$ of .002.

Clearly the irregularization has substantially slowed the movement of the shock and thereby induced a delta-like spike to conserve mass.  The extent to which this is a bona fide delta shock is a delicate question (we will thus cautiously refer to this feature as a \textit{delta spike}).  We mention three relevant considerations:

\begin{enumerate}
\item To begin with, the effect is inherently non-smooth, as for smooth solutions of eq.~(\ref{modHopf2}), we can prove the following maximum principle:
\begin{prop}
A solution $u$ to eq.~(\ref{modHopf2}) on a compact domain satisfies
\begin{equation}
max\{u(\cdot,t)\} = max\{u(\cdot,0)\}
\end{equation}
for all times $t>0$ while the solution remains smooth.
\end{prop}

We motivate this as follows: If we look at a curve $y(t)$ such that $u(\cdot,t)$ attains a maximum at $y(t)$, then we find
\begin{eqnarray}
\partial_t u(y(t),t) & = & u_x \dot{y} + u_t \nonumber \\
& = & u_x \dot{y} - \left(\frac{\frac{1}{2}u^2}{1+\epsilon \: u_x^2}\right)_x \nonumber \\
& = & u_x \dot{y} - \left(\frac{u_x u }{1+\epsilon \: u_x^2} - \frac{\epsilon u_x u_{xx} u^2}{\left(1+\epsilon \: u_x^2\right)^2}\right) \nonumber \\
& = & 0
\label{maxprinciple}
\end{eqnarray}
Where the last expression vanishes because $u_x=0$ at the maximum $y(t)$.  This shows that the maximum value of $u$ is unchanging in time.  (This argument can be made more rigorous, avoiding the assumption that the maximum can be parameterized by a differentiable curve $y(t)$, but we will not do so here.)  Note that a similar result holds also in the case of eq.~(\ref{modHopf3}) with a 2nd order (viscous) regularization, with the modification that the maximum of $u$ is then decreasing in time.  

\item Examining the delta spike closely reveals an internal structure.  To provide some perspective, we note first that eq.~(\ref{modHopf2}) admits two relevant exact solutions, one of the form $a(t) + b(t) x$ (i.e. a straight line, moving in time), and the other a stationary solution of the form
\begin{equation}
\alpha \cosh\left(\frac{x-\beta}{\alpha \sqrt{\epsilon}}\right).
\label{coshProfile}
\end{equation}

Experimentally, we find that the shape of the spike is well approximated by a profile of the form (\ref{coshProfile}), as shown in fig.~(\ref{coshFig}).  More precisely, the form of the delta spike consists of two partial coshes, one positive and one negative, separated by an abrupt discontinuity.  Away from the defect, the solution to eq.~(\ref{modHopf2}) approaches the aforementioned linear solution asymptotically in time.  This eventual shape is analogous to the asymptotic ``N-wave'' profile of solutions to the regular Hopf equation.  As discussed in the appendix, it is possible to derive analytic expressions for the parameters $a(t)$, $b(t)$, $\alpha$, $\beta$ describing this asymptotic form, and we find good agreement between simulation and these analytic values. 

\begin{figure}[h]
	\includegraphics[width=1.1\columnwidth,center]{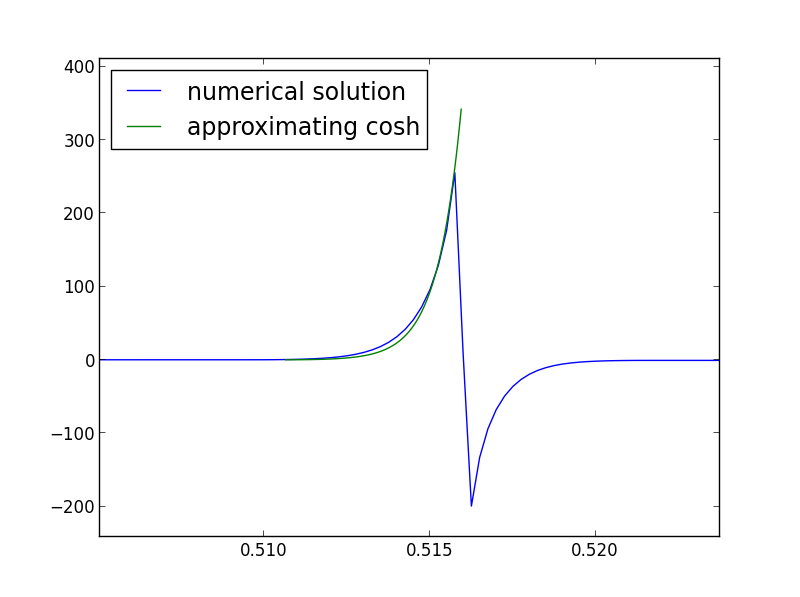}
	\caption{{\bf Cosh profile of delta spike.} Comparison of a delta spike with a cosh profile of form (\ref{coshProfile}).  Note that the cosh profile shown here is not a best fit, but rather an analytically derived approximation, as discussed in the appendix.  The delta spike comes from a simulation with parameters $\epsilon=2.5e-6,\;dx=2.5e-4$, and initial data $u(x,0)=\sin(2\pi x)+.1$. Again, the positive and negative spikes must be added together to give the net
weight of the delta shock.
	}
	\label{coshFig}
\end{figure}

An important feature of the cosh profile (\ref{coshProfile}) is that its width scales like $\sqrt{\epsilon}$.  To conserve mass, the height of the shock must therefore scale like $1/\sqrt{\epsilon}$.  This suggests that, although the delta spike is not truly a delta \textit{shock} for non-zero values of $\epsilon$ (indeed, the cosh profile (\ref{coshProfile}) certainly does not qualify as a delta shock), as $\epsilon \rightarrow 0$ the spike does indeed converge to a delta function.  Experimentally, this is what we seem to find (conditional upon some finessing of the numerical methods), as discussed in the next point. 

\item Establishing numerical convergence of our simulations is tricky for two reasons: Firstly, since we are looking for convergence to a delta shock, we cannot hope to have convergence in e.g. the $L_2$ sense.  Secondly, the equation (\ref{modHopf2}) is highly unstable (not surprisingly, since the irregularization is designed specifically to produce a delta shock, which may be viewed as a sort of instability).  After trying various finite-differencing schemes, we found that the best results came from a type of conservative upwind differencing, wherein $u$ is differenced upwind but $u_x$ is not.  Explicitly, the semi-discrete formulation of eq.~(\ref{modHopf1}) is 
\begin{gather}
\partial_t u_i + \Delta F_{i-1} = 0 \nonumber \\
\nonumber \\
F_i = \left\{ 
\begin{array}{l l}
u_i^2 f\left(\epsilon,\frac{\Delta u_i}{dx}\right) & \quad (u_i+u_{i+1}) > 0\\
u_{i+1}^2 f\left(\epsilon, \frac{\Delta u_i}{dx}\right) & \quad (u_i+u_{i+1}) < 0 
\end{array} 
\right.
\end{gather}
where $\Delta$ is the discrete forward difference operator ($\Delta a_i=a_{i+1}-a_i$).  All figures presented herein were produced using this discretization.  (Remark: for the regularized eq.~(\ref{modHopf3}) we performed an operator splitting, evaluating the regularization term in Fourier space for stability.)  Despite being the most successful numerical method we found, the stability properties of this discretization are not entirely satisfactory.  In particular, we found that when the grid size became substantially smaller than the width of the delta spike (which is on the order of $\sqrt{\epsilon}$), the delta spike would become unstable and break apart into multiple, smaller delta spikes.  To overcome this difficulty, we let the mesh size $dx$ and the width of the delta $\sqrt{\epsilon}$ go to zero together, keeping $dx/\sqrt{\epsilon}$ fixed (recall that it was our intention to let $\epsilon\rightarrow 0$ anyhow, so as to recover a solution to the original Hopf equation).  This technique appears to work, as seen in fig.~(\ref{deltaCvrg}):  The delta spike gets narrower and higher, while remaining in essentially the same place, as $dx,\epsilon \rightarrow 0$ in this manner.

\begin{figure}[h!]
	\includegraphics[width=1.1\columnwidth]{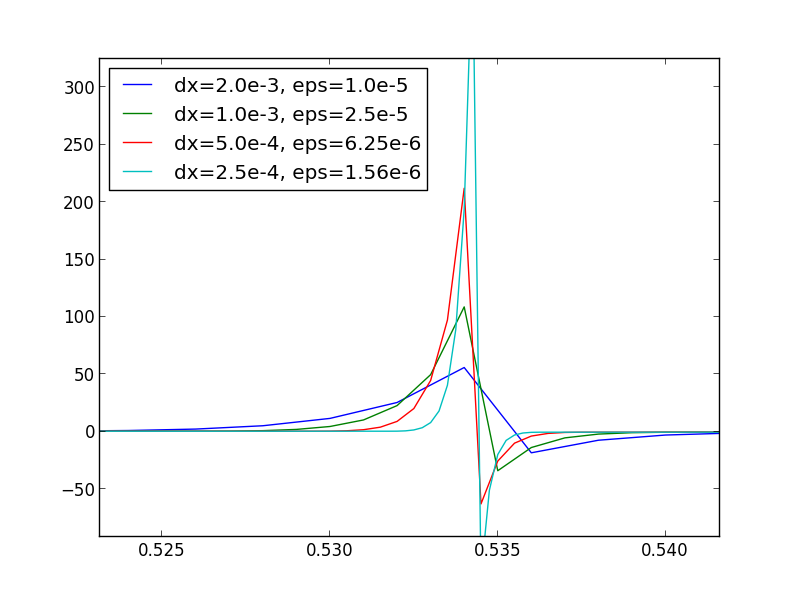}
	\caption{(Color online) Close-up of the delta defects as $dx,\epsilon\rightarrow 0$, $\sqrt{\epsilon}\sim dx$.
	}
	\label{deltaCvrg}
\end{figure}

\end{enumerate}

Although our simulations seem to establish convergence to a delta shock, it is a very peculiar convergence wherein the strength $\epsilon$ of the irregularization must vanish together with the discretization length $dx$.  The fact that the width of the delta spike must be kept close to $dx$ suggests that there may be some sort of resonance between the irregularization of eq.~(\ref{modHopf2}) and the grid which drives the formation of the spike (especially in light of the maximum principle mentioned previously).  Nonetheless, we can still hope to produce a meaningful weak solution for the Hopf equation in this way.  We will not establish anything more rigorous at this time, though. 


\section{Physical motivation: Delta shocks and dislocation walls}
\label{sec:Physics}

Our interest in delta shocks is motivated by a tangible physical question.
Why do dislocations in crystals form walls~\cite{LimkumnerdSPRL,ChenPRL}? 
Briefly, a crystal is a regular
array of atoms. A dislocation line is a flaw in that array, such as the
edge dislocation formed by the boundary of an extra plane of atoms, or the 
screw dislocation forming the central line where planes form
a `spiral staircase'~\cite{HirthLothe}. These dislocations move to mediate
plastic deformation when the crystal is bent, and form tangles that
organize into wall-like structures called {\em cell walls}~\cite{Hughes1998}.
In a continuum theory of dislocation dynamics, such walls must be described
as delta shocks: their density scales as the inverse square of the lattice
constant, which vanishes in the continuum limit.

One class of continuum dislocation dynamics
theories~\cite{Roy2005,LimkumnerdSPRL} do form such delta shocks. In these
theories, the dislocations move with a common velocity $\cal{F}$. If one
describes an incipient wall in the $yz$ plane with a dislocation density
that depends only on $x$, the dependence of this velocity on the dislocation
density simplifies~\cite{Limkumnerd2006} to the Hopf equation:
\begin{equation}
\mathcal{F}_t + \left(\frac{1}{2} \mathcal{F}^2\right)_x = 0.
\label{PKF}
\end{equation}
When the Hopf equation forms a step-like shock, the entrained dislocations
pile up into a delta shock, forming a wall.

The central physical question is how that wall should evolve after it forms.
If we regularize the Hopf equation~(\ref{PKF}) into Burger's equation (i.e. viscous regularization),
the dislocation wall moves along with the Rankine-Hugoniot velocity.
This produces a self-consistent, sensible model~\cite{ChenPRL,Chen2013},
but one that is unsatisfying in two regards. First, the velocity of
dislocation walls in crystals is determined by the microstructure of the
walls and not by continuum properties. Secondly, generalizations of 
this continuum theory have differing velocities for the different components
of the dislocation density -- smearing the resulting walls.

On the other hand, \textit{physically} at some junctions where two dislocations intersect, the 
point of intersection can be pinned in place (so-called {\em sessile}
dislocation junctions~\cite{HirthLothe}). Dislocation walls can also
be pinned by impurities that segregate to the boundary. This provides the
motivation for our study in section~(\ref{sec:Stationarity})
of the Hopf equation~(\ref{PKF}) with stationary shocks. 

We should note that the behavior of both our model in two and three dimensions
and experimental dislocation systems under stress is more complex than
simple wall singularities. There they produce complex cellular structures,
which in our models and some experiments form self-similar, fractal
morphologies. Our theories in higher dimensions show clear
analogies~\cite{WoosongTurbulence} to the behavior of the Euler
equation (the inviscid limit of Navier-Stokes). Nonetheless, these
structures are complex, ramified wall-like entities, whose dynamics
should be controlled by physics on the microscopic atomic scale, not by the 
continuum laws.

\section{Conclusions}
We have seen that delta shocks arise naturally in a variety of systems of conservation laws, and their presence can often be traced to an incompatibility of the flow of the conserved quantity and the motion of defects.  Although most of the existing examples of delta shocks do not have direct physical interpretation, the case of cell walls in crystals does appear to furnish such an example.  

A recurring theme in our presentation has been the importance of developing numerical techniques capable of handling defects with more exotic behavior than traditional, viscous shocks.  Ideally, one would like to be able to specify whatever defect behavior is believed appropriate for a given system and have numerical techniques which respect that behavior.  In lieu of such a very general approach to simulating conservation laws, we must instead focus on specific classes of defect behavior which we believe are important.  To this end, we propose further study of the stationarity condition introduced above.  

We have concocted an ``irregularization'' (\ref{modHopf1}-\ref{modHopf3}) of the Hopf equation which appears numerically to exhibit delta shocks, suggesting a possible avenue to development of the aforementioned numerical methods for delta shocks.  Though we are far from showing that these apparent deltas are \textit{bona fide}, we nonetheless think that the irregularization and the defects it yields are interesting in their own right and merit continued study.


\section{Appendix}
We present here some further numerical and analytical observations (mostly without details or proofs) for the irregularization
\begin{equation}
u_t + \left(\frac{\frac{1}{2}u^2}{1+\epsilon \: u_x^2}\right)_x = 0.
\label{rairreg}
\end{equation}
In particular, we sketch a systematic way to approximate the asymptotic form of solutions to this equation. 

Three useful characteristics of the regular Hopf equation can be shown to also hold for smooth solutions of this more general equation: 
\begin{enumerate}
\item The maximum principle mentioned above.
\item Zero-crossings of a solution to eq. (\ref{rairreg}) remain fixed in time.
\item The total mass between any two zero-crossings is conserved.
\end{enumerate}
The second and third items above allows us to predict (analytically, approximately) the asymptotic form of solutions to eq.~(\ref{rairreg}).  We sketch this analysis as follows:  Shocks tend to form on downwards-going zero crossings, and away from shocks the solution becomes linear at large times.  The overall profile of the solution thus is linear everywhere except at the shocks, where it is approximated by the cosh profile mentioned above.  The stationarity of \textit{upwards}-going zero crossings allows for analytical evaluation of the slope there, which asymptotically must equal the slope of the entire linear region of the solution.  Explicitly, the slope $u_x$ at a zero crossing $x_{zc}$ satisfies the ODE
\begin{eqnarray}
0 & = & \partial_t u_x(x_{zc},t) + \frac{1}{2}\left(\frac{u(x_{zc},t)^2}{1+\epsilon u_x(x_{zc},t)^2}\right)_{xx} \nonumber \\
& = & \partial_t u_x(x_{zc},t) + \left(\frac{u_x(x_{zc},t)^2}{1+\epsilon u_x(x_{zc},t)^2}\right)
\label{ODE}
\end{eqnarray}
which is obtained by differentiating eq.~(\ref{rairreg}) w.r.t. $x$ and noting that all terms with a factor of $u$ vanish.  This ODE (\ref{ODE}) can be solved analytically, providing the asymptotic form of the linear part of the solution.

The third item above allows us to evaluate the mass of the delta, which is just the difference between the total initial mass between two zero-crossings and the eventual mass of the linear solution on the same region.  Examining the cosh profile of the delta mentioned above, 
\begin{equation*}
\alpha \cosh\left(\frac{x-\beta}{\alpha \sqrt{\epsilon}}\right),
\end{equation*}
we see that $\alpha$ must equal the height at which the delta begins, which can be determined from the linear solution.  Thus $\alpha$ is determined, and then $\beta$ can be determined from the known mass of the delta.  

Stringing together the above arguments allows a piecewise description of the solution across the entire domain.  This analytic, approximate solution is found to agree reasonably well with simulation---for instance, the cosh profile shown in fig.~(\ref{coshFig}) uses parameters calculated in this way.  

\section{Acknowledgement}
We thank Randall LeVeque for helpful comments and guidance on hyperbolic conservation laws and numerical methods.  We are especially indebted to Alexander Vladimirsky 
for his
exceptionally insightful suggestions for the production of this paper.  This work was supported by US DOE/BES grant DE-FG02-07ER46393.


\begin{thebibliography}{99}
\bibitem{Lax1973}
	P. Lax,
	\emph{Hyperbolic Systems of Conservation laws and the Mathematical Theory of Shock Waves},
	Society for Industrial and Applied Mathematics,
	1973

\bibitem{Lax1957}
	P. Lax,
	\emph{Hyperbolic Systems of Conservation Laws II},
	Communications on Pure and Applied Mathematics,
	1957

\bibitem{LeFloch2002}
	P. G. Le Floch,
	\emph{Hyperbolic Systems of Conservation Laws: The Theory of Classical and Nonclassical Shock Waves},
	Birkhäuser Verlag AG, c/o Springer GmbH \& Co,
	2002

\bibitem{Evans}
	L. Evans,
	\emph{Partial Differential Equations, second ed.},
	American Mathematical Society,
	1998

\bibitem{Tan1994}
	{ D. C. Tan, T. Zhang, T Chang, and Y. X. Zheng},
	{Delta-Shock Waves as Limits of Vanishing Viscosity for Hyperbolic Systems of Conservation Laws},
	{Journal of Differential Equations},
	{1994}

\bibitem{Colombeau1984}
	J. F. Colombeau,
	\emph{New Generalized Functions and Multiplication of Distributions },
	Elsevier Science Publishers B.V., 
	1984

\bibitem{Keyfitz1989}
	B. L. Keyfitz and H. C. Kranzer
	\emph{A Viscosity Approximation to a System of Conservation Laws with No Classical Riemann Solution}
	"Nonlinear Hyperbolic Problems," Lecture Notes in Mathematics, 
	Springer-Verlag,
	1989

\bibitem{Keyfitz1999}
	B. L. Keyfitz,
	\emph{Conservation Laws, Delta Shocks and Singular Shocks},
	Nonlinear Theory of Generalized Functions,
	1999

\bibitem{Panov2006}
	E. Yu. Panov and V. M. Shelkovich,
	\emph{$\delta'$-Shock waves as a new type of solutions to systems of conservation laws},
	Journal of Differential Equations,
	2006

\bibitem{Meppelink2009}
	R. Meppelink, S. B. Koller, J. M. Vogels, P. van der Straten, E. D. van Ooijen, N. R. Heckenberg, H. Rubinsztein-Dunlop, S. A. Haine, and M. J. Davis,
	\emph{Observation of shock waves in a large Bose-Einstein condensate},
	Physical Review A,
	2009

\bibitem{Taylor1970}
	R. J. Taylor, D. R. Baker, and H. Ikezi,
	\emph{Observation of collisionless electrostatic shocks},
	Physical Review Letters,
	1970


\bibitem{Kurganov2002}
	{Kurganov, A.  and Noelle, S.  and Petrova, G.},
	{Semidiscrete central-upwind schemes for hyperbolic conservation laws and Hamilton--Jacobi equations},
	{SIAM J. Sci. Comput.},
	{2001}

\bibitem{HirthLothe}
	J. P. Hirth and J. Lothe,
	\emph{Theory of Dislocations, 2nd ed.},
	Krieger Publishing Company,
	1982

\bibitem{Limkumnerd2006}
	S. Limkumnerd and J.P. Sethna,
	\emph{Mesoscale theory of grains and cells: crystal plasticity and coarsening},
	Physical Review Letters,
	2006

\bibitem{Roy2005}
	A. Roy and A. Acharya,
	\emph{Finite element approximation of field dislocation mechanics},
	J. Mech. Phys. Solids, 
	2005

\bibitem{Hughes1998}
	D. A. Hughes, D. C. Chrzan, Q. Liu, and N. Hansen,
	\emph{Scaling of Misorientation Angle Distributions},
	Physical Review Letters,
	1998

\bibitem{Choi2013}
	W. Choi,
	\emph{The Physics of Singular Dislocation Structures in Continuum Dislocation Dynamics},
	Thesis, Cornell University,
	2013


\bibitem{Chen2013}
	Y. S. Chen, W. Choi, M. Bierbaum, J. P. Sethna, and S. Papanikolaou,
	\emph{Scaling theory of continuum dislocation dynamics in three dimensions: Self-organized fractal pattern formation},
	International Journal of Plasticity,
	2013

\bibitem{LimkumnerdThesis}
	S. Limkumnerd,
	\emph{Mesoscale Theory of Grains and Cells: Polycrystals \& Plasticity}
	Thesis, Cornell University
	2007

\bibitem{LimkumnerdSPRL}
	S. Limkumnerd and J. P. Sethna, 
	\emph{Mesoscale Theory of Grains and Cells: Crystal Plasticity and
	Coarsening}, Phys. Rev. Lett. {\bf 96}, 095503 (2006).

\bibitem{ChenPRL}
	Yong S. Chen, Woosong Choi, S. Papanikolaou, and J. P. Sethna,
	\emph{Bending Crystals: Emergence of Fractal Dislocation Structures},
	Phys. Rev. Lett. {\bf 105}, 105501 (2010).

\bibitem{WoosongTurbulence} 
        Woosong Choi, Yong S. Chen, Stefanos Papanikolaou, and James P. Sethna,
        \emph{Is dislocation flow turbulent in deformed crystals?},
	Computing in Science and Engineering,
        {\bf 14}, 33 (2012).


\end{thebibliography}
\end{document}